\documentclass{amsart}
\usepackage{amssymb,amscd}  

\newtheorem{theorem}{Theorem}[section]
\newtheorem{lemma}[theorem]{Lemma}
\newtheorem{corollary}[theorem]{Corollary}
\newtheorem{conjecture}[theorem]{Conjecture}
\newtheorem{proposition}[theorem]{Proposition}

\theoremstyle{definition}
\newtheorem{definition}[theorem]{Definition}

\newtheorem{remark}[theorem]{Remark} 
\numberwithin{equation}{section}

\newcommand\C{\mathbf{C}} 

\newcommand\Q{\mathbf{Q}}
\newcommand\Z{\mathbf{Z}}
\newcommand\F{\mathbf{F}}
\newcommand\twomatr[4]{\begin{pmatrix}#1&#2\\ #3&#4 \end{pmatrix}}
\newcommand\stwomatr[4]{\left(\begin{smallmatrix}#1&#2\\  
                               #3&#4 \end{smallmatrix}\right)}
\newcommand\twovec[2]{\begin{pmatrix}#1\\ #2 \end{pmatrix}}
\newcommand\stwovec[2]{\left(\begin{smallmatrix}#1\\  
                               #2 \end{smallmatrix}\right)}

\newcommand\tensor{\otimes}
\newcommand\isomorphic{\cong}

\newcommand\Directsum{\bigoplus}
\newcommand\directsum{\oplus}
\newcommand\union{\cup}

\newcommand\Intersect{\bigcap}
\newcommand\abs[1]{{\left|#1\right|}}

\newcommand\Projective{{\bf P}} 
\newcommand\modforms{\mathcal{M}}

\newcommand\Half{\mathcal{H}} 
\newcommand\LL{\mathcal{L}} 
\newcommand\thet{\theta}    
\newcommand\phee{\varphi}   
\newcommand\morph{F}        
\newcommand\Gammapm{\Gamma_{\pm}} 
\newcommand\Xpm{X_{\pm}}    

\begin{document}

\title{On the maps from $X(4p)$ to $X(4)$}

\author{Samar Jaafar}
\address{Department of Mathematics, Rutgers University, Hill Center, 110
  Frelinghuysen Road, Piscataway, NJ 08854, USA}
\email{jaafar@math.rutgers.edu}

\author{Kamal Khuri-Makdisi}
\address{Corresponding author: Mathematics Department and Center for
Advanced Mathematical Sciences, American University of Beirut, Bliss
Street, Beirut, Lebanon}
\email{kmakdisi@aub.edu.lb}

\thanks{February 28, 2008}

\begin{abstract}
We study pullbacks of modular forms of weight $1$ from the modular curve
$X(4)$ to the modular curve $X(4p)$, where $p$ is an odd prime.  We find
the extent to which such modular forms separate points on $X(4p)$.
Our main result is that these modular forms give rise to a morphism $F$
from the quotient of $X(4p)$ by a certain involution $\iota$ to projective
space, such that $F$ is a projective embedding of $X(4p)/\iota$ away
from the cusps.  We also report on computer calculations regarding products
of such modular forms, going up to weight $4$ for $p \leq 13$, and up to
weight $3$ for $p \leq 23$, and make a conjecture about these products and
about the nature of the singularities at the cusps of the image
$F(X(4p)/\iota)$.
\end{abstract}

\maketitle

\section{Introduction}
\label{section1}

For a congruence subgroup $\Gamma \subset \Gamma(1) = SL(2,\Z)$, we denote
by $\modforms_k(\Gamma)$ the space of modular forms
of weight $k \geq 0$ with respect to $\Gamma$, which we view as
holomorphic functions on the upper half plane $\Half$.  We denote by
$X(\Gamma)$ the modular curve obtained by completing the quotient
$\Gamma\backslash\Half$ at the cusps.  For the principal congruence
subgroup $\Gamma(N)$ consisting of matrices congruent to
$\stwomatr{1}{0}{0}{1}$ modulo $N$, we abbreviate $X(N)$ for the modular
curve $X(\Gamma(N))$.

In this article, we study a linear series on the modular curve $X(4p)$
corresponding to a subspace $V \subset \modforms_1(\Gamma(4p))$, for $p$ an
odd prime.  The space $V$ is spanned by forms of the form $f \mid \alpha$,
where $f \in \modforms_1(\Gamma(4))$, and $\alpha$ runs over the set of
integral matrices with determinant $1$ or $p$.  We can view the different
$\alpha$ as maps from $X(4p)$ to $X(4)$, induced by $z \mapsto \alpha z \in
\Half$; thus $V$ is spanned by the pullbacks of forms from $X(4)$.  Now our
incomplete linear series $V$ gives rise to a morphism $\morph$ from $X(4p)$
to projective space, that reflects all the maps from $X(4p)$ to $X(4)$ as
well as some additional information.  This morphism factors through an
involution $\iota$ of $X(4p)$, giving rise to a morphism (which we still
call $\morph$) from $\Xpm(4p) = X(4p)/\iota$ to projective space.  Our main
results state that $\morph$ gives a birational projective embedding of
$\Xpm(4p)$ that is an isomorphism away from the cusps
(Theorem~\ref{theorem3.2}), and that the effect of $\morph$ on the cusps is
to partially identify the cusps in a precise manner
(Theorem~\ref{theorem3.5}).  This identification of the cusps gives rise to
singularities on the image curve $\morph(\Xpm(4p))$, and we formulate a
conjecture (Conjecture~\ref{conjecture1.1})
about the nature of these singularities in terms of the
algebra of modular forms generated by elements of $V$, under multiplication
of forms.  Since $V \subset \modforms_1(\Gamma(4p))$, a product of $k$
elements of $V$ is naturally a weight $k$ form in
$\modforms_k(\Gamma(4p))$.  Such weight $k$ forms correspond to the
pullback via $\morph$ of degree $k$ polynomials from projective space to
$\Xpm(4p)$. 
The fact that the morphism $\morph$ identifies cusps of $\Xpm(4p)$ in a
certain way means that modular forms in $V$, and hence also their products
in $\modforms_k(\Gamma(4p))$, assume equal values at cusps that have been
identified; this imposes linear conditions on the forms that we obtain in
$\modforms_k(\Gamma(4p))$, above and beyond the conditions caused by the
involution $\iota$.  We then conjecture that the singularities of
$\morph(\Xpm(4p))$ at the images of cusps involve only first-order contact
at the cusps that are identified with each other; if there were
higher-order contact, this would imply further linear relations between
higher-order coefficients in the Taylor (i.e., $q$-) expansions of our
forms at such cusps for all weights $k$.  In algebraic terms, our
conjecture can be stated as follows:

\begin{conjecture}
\label{conjecture1.1}
For $k \geq 4$, the degree $k$ monomials in elements of $V$ span a
subspace $W_k \subset \modforms_k(\Gamma(4p))$ whose dimension is
\begin{equation}
\label{equation1.1}
\dim W_k = 
(1/2) \dim \modforms_k(\Gamma(4p)) - 3(p+1)(p-3)
 = (k-1)p^3 - 3p^2/2 + (7-k)p + 15/2.
\end{equation}
The factor $1/2$ corresponds to invariance under the involution $\iota$,
and the ``missing'' $3(p+1)(p-3)$ dimensions correspond to the
identifications among cusps due to $\morph$.
\end{conjecture}
Once the dimension is attained for a given $k$, we expect it to be attained
for all larger $k$.  We have numerically computed the dimensions of the
$W_k$ for $k \leq 4$ when $p \leq 13$, and for $k \leq 3$ for the
remaining $p \leq 23$.
The above computations required more computer time than one might expect:
e.g., for the case $p = 23$, it took 16 days for a desktop computer to go
up to $k=3$, using our code which involved C++ code with the NTL
library~\cite{NTL} for the linear algebra, and some MAGMA~\cite{MAGMA} code
to compute a Gr\"obner 
basis at one point.  The last and lengthiest step in this calculation,
which was to find $\dim W_3$, involved linear algebra on a matrix of size
$50,000 \times 30,000$, with integral entries.  Thankfully, the matrix
could be recast as a block diagonal matrix with 92 blocks, each of size
approximately $560 \times 325$, whose ranks we needed to find; these were
dense integral matrices that started out with entries of about $7$ digits,
but whose entries grew considerably during the linear algebra due to
coefficient explosion.  In Section~\ref{section4}, we sketch how we did the
computation, and present the results that suggested
Conjecture~\ref{conjecture1.1}.

The second named author would like to thank the first named author for
carrying out the computer calculations described in Section~\ref{section4}.
The second named author also gratefully acknowledges support from the
University Research Board, at the American University of Beirut, and from
the Lebanese National Council for Scientific Research, for supporting this
work through the grants ``Equations for modular and Shimura curves.''

\section{The space $V$ of modular forms}
\label{section2}

To motivate our choice of $V$, we first discuss a principle that allows one
to obtain projective models of curves, in a form that is suitable for
computation.  If $X$ is a smooth projective curve of genus $g$ over $\C$,
and $\LL$ is a line bundle on $X$ with $\deg \LL \geq 2g+2$, then the
complete linear series on $X$ arising from $H^0(X,\LL)$ gives a projective
embedding of $X$ that is projectively normal, and moreover the ideal
defining the image of $X$ in projective space is generated by
quadrics~\cite{SaintDonatEquations}.  Thus the pair $X$ and $\LL$ are
completely described by a knowledge both of the spaces $H^0(X,\LL)$ and
$H^0(X,\LL^{\tensor 2})$, and of the multiplication map  
$\mu: H^0(X,\LL) \tensor H^0(X,\LL) \to H^0(X,\LL^{\tensor 2})$.  
This model of $X$ allows for efficient algorithmic implementation
of the group operations in the Jacobian of $X$, as
in~\cite{KKMlinear, KKMasymptotic}.
When $X = X(\Gamma)$ is a modular curve, 
we can choose $\LL$ so that $H^0(X,\LL) = \modforms_k(\Gamma)$ and
$H^0(X,\LL^{\tensor 2}) = \modforms_{2k}(\Gamma)$ for a suitable weight
$k$.  If the modular curve $X(\Gamma)$ has at least $4$ cusps (as is the
case for $X(N)$ for $N \geq 3$), it is immediate that $k = 2$ corresponds
to a line bundle $\LL$ with $\deg \LL \geq 2g+2$; thus we can use any $k
\geq 2$ to produce the model of our modular curve.  The multiplication map
$\mu$ is then multiplication of modular forms, which can be carried out
using $q$-expansions, where one can work only with $q$-expansions to a
fixed order $O(q^n)$, using a sufficient number of coefficients to be able
to recognize a form in the finite-dimensional space
$\modforms_{2k}(\Gamma)$.  Instead of working with $q$-expansions, we can
work with the values of the modular forms at sufficiently many points of
$\Half$, as in ongoing investigation of the second named author; thus this
approach may eventually lead to explicit models of Shimura curves, even in
the absence of cusps and of $q$-expansions --- provided that we can get a
computational handle on the entire space $\modforms_k(\Gamma)$.  This poses
no problem in the classical case, where we can compute $q$-expansions for a
basis of the space $\modforms_k(\Gamma)$ for $k\geq 2$ using modular
symbols.

In the setting of Shimura curves, however, it appears to be more difficult
to compute the entire space $\modforms_k(\Gamma)$ for $\Gamma$ an
arithmetic group associated to a quaternion algebra, especially as the
``level'' of $\Gamma$ increases.  It may be possible in general to find
just a few modular forms of fixed small level, and then to use these to
produce more forms of higher varying level.  An example of this in the
elliptic modular case is to start with the Eisenstein series $E_4(z)$ and
$E_6(z)$ on $\Gamma(1)$, and then to consider forms such as
$E_4((az+b)/d)$ and $E_6((az+b)/d)$ with $ad | N$, which belong to
$\modforms_4(\Gamma(N))$ and $\modforms_6(\Gamma(N))$.  Of course, the
subspaces that such forms span in $\modforms_4(\Gamma(N))$ and in
$\modforms_6(\Gamma(N))$ fall far short of the entire spaces of modular
forms on $\Gamma(N)$ in weights $4$ and~$6$: for instance, we do not obtain
any cusp forms in these weights.  On the other hand, by taking
\emph{products} of such forms, we can obtain a much larger set of modular
forms --- in general, taking products of modular forms mixes up the
different automorphic representations to such an extent that the Hecke
orbit of a product often generates the full space of modular forms in a
given weight and level.

In this article, we study a specific case of this way of generating
elliptic modular forms by taking products, with the goal of observing
phenomena that will be of interest in future research related to Shimura
curves.  Even in this setting, we have come up against some conceptual and
computational challenges.  We shall not work with $E_4$ and $E_6$, but will
prefer to work with the slightly higher ``base'' level of $\Gamma(4)$, so
as to be able to use modular forms of weight~$1$ as our simple forms.  The
smaller weight means that the coefficients of the $q$-expansions are
correspondingly smaller, which has helped us considerably with our computer
calculations involving products of modular forms.  It is also relevant that
$\Gamma(4)$ does not contain any elliptic elements, so we avoid problems
arising from excess vanishing of modular forms at elliptic points.

All the forms that we shall encounter can be expressed
in terms of the following modular forms on $\Gamma(4)$, where as usual 
$q = \exp(2\pi i z)$:
\begin{equation}
\label{equation2.1}
\begin{split}
&\thet(z) = \sum_{n \in \Z} q^{n^2} 
         = \sum_{\text{even numbers } k \in \Z} q^{k^2/4}
         = 1 + 2q + 2q^4 + \cdots,\\ 
&\phee(z) = \thet(z) - \thet((z+2)/4) 
         = \sum_{\text{odd numbers } k \in \Z} q^{k^2/4}
         = 2q^{1/4} + 2q^{9/4} + \cdots, \\
&M(z) = \thet(z)^2,
\qquad
N(z) = 2 \thet(z) \phee(z),
\qquad
P(z) = \phee(z)^2.
\end{split}
\end{equation}
Here $\thet$ and $\phee$ have weight $1/2$, and $M$, $N$, and $P$ have
weight~$1$.
\begin{proposition}
\label{proposition2.1}
The set $\{M, N, P\}$ is a basis for the space $\modforms_1(\Gamma(4))$,
and the graded ring $\directsum_{k \geq 0} \modforms_k(\Gamma(4))$ is
generated as a $\C$-algebra by $M$, $N$, and $P$, modulo the relation 
$N^2 - 4MP = 0$.
\end{proposition}
\begin{proof}
Consider the subring of $\directsum_{k \geq 0} \modforms_k(\Gamma(4))$
generated by $M$, $N$, and $P$.  The degree $k$ part of this subring
is spanned by the monomials $\thet^{2k-j} \phee^j = 2^j q^{j/4} + \cdots$
for $0 \leq j \leq 2k$.  These are plainly linearly independent, so the
degree $k$ part of our subring has dimension $2k+1$.  On the other hand,
the curve $X(4)$ is of genus $0$, with $6$ cusps and without any elliptic
points.  This implies that $\dim \modforms_k(\Gamma(4)) = 2k+1$ for all
$k \geq 0$ (see sections 1.6 and~2.6 of~\cite{Shimura}, including equation
(2.6.9) for the case $k=1$).  Thus our subring is indeed the whole
graded ring of modular forms on $\Gamma(4)$.  By considering the Hilbert
series of this graded ring, we see that there can be no relation between
$M$, $N$, and $P$ besides the obvious one $N^2 - 4MP = 0$.
\end{proof}

\begin{remark}
\label{remark2.2}
Proposition~\ref{proposition2.1} says that $X(4)$ is isomorphic to the
projective line $\Projective^1$, and that the projective embedding $X(4)
\to \Projective^2$ arising from $\modforms_1(\Gamma(4))$ comes from the
line bundle $\LL = \mathcal{O}_{\Projective^1}(2)$, realizing $X(4)$ as a
projective plane conic.   We have not tried to 
work directly with the weight $1/2$ forms $\thet$ and $\phee$ because it is
harder with half-integral weight to keep track of factors of automorphy
when one works with forms such as $\thet((z+b)/p)$ or $\phee(pz)$ of level
$4p$.  We note incidentally that $\mathcal{S}_1(\Gamma(4)) = 0$, which
means that the forms $M$, $N$, and $P$ are Eisenstein series.
\end{remark}

We now let $p \geq 3$ be a prime, and propose to study forms of the form
$f|\alpha$, with $f \in \modforms_1(\Gamma(4))$, for which 
$f|\alpha \in \modforms_1(\Gamma(4p))$.  Since we can multiply $\alpha$ by
a scalar, we shall always assume that $\alpha$ is a primitive matrix, i.e.,
an integral matrix with relatively prime entries.  We thus define our basic
vector space $V$ of weight~$1$ forms on $\Gamma(4p)$ as follows:

\begin{definition}
\label{definition2.3}
We define $V \subset \modforms_1(\Gamma(4p))$ to be the span of all modular
forms of the form $(f | \alpha)(z) = f(\alpha z) j(\alpha, z)^{-1}$, 
where $f \in \modforms_1(\Gamma(4))$, $\alpha = \stwomatr{a}{b}{c}{d}$ is a
primitive matrix such that $\alpha \Gamma(4p) \alpha^{-1} \subset
\Gamma(4)$ and $\det \alpha > 0$, and $j(\alpha, z) = cz + d$.  By
Proposition~\ref{proposition2.1}, we can always assume that
$f \in \{M, N, P\}$. 
\end{definition}
Note that a primitive $\alpha$ as in the above definition gives rise to a
map from $X(4p)$ to $X(4)$, induced from $z \mapsto \alpha z$ on $\Half$.
Also note that we will only need the slash operator $f \mapsto f|\alpha$
for forms of weight~$1$, so we do not need the more general factor
$j(\alpha, z)^{-k}$ for weight~$k$.

\begin{lemma}
\label{lemma2.4}
The matrices $\alpha$ of Definition~\ref{definition2.3} are precisely those
primitive matrices with determinant $1$ or $p$.  Moreover, define the finite
set
\begin{equation}
\label{equation2.2}
S = \left\{ \twomatr{1}{0}{0}{1}, 
\twomatr{1}{0}{0}{p}, \twomatr{1}{1}{0}{p}, \twomatr{1}{2}{0}{p},
\dots, \twomatr{1}{p-1}{0}{p} \right\}.
\end{equation}
Then we have
\begin{equation}
\label{equation2.3}
\begin{split}
V  &= \mathrm{span}\bigl\{ f|\alpha \mid f \in \{M, N, P\} \text{ and }
\alpha \in S\bigr\} \\
   &= \mathrm{span} \{ f(z), f((z+b)/p)
            \mid  f \in \{M, N, P\} \text{ and } 0 \leq b < p\}.\\
\end{split}
\end{equation}
Hence $\dim V \leq 3(p+1)$.
\end{lemma}
\begin{proof}
Observe first that $\Gamma(4)$ and $\Gamma(4p)$ are normal subgroups of
$\Gamma(1)$.  It follows that
the condition $\alpha \Gamma(4p) \alpha^{-1} \subset \Gamma(4)$
depends only on the double coset $\Gamma(1) \alpha \Gamma(1)$.  Now this
double coset is equal to $\Gamma(1) \stwomatr{\ell}{0}{0}{1} \Gamma(1)$,
where $\ell = \det \alpha$, and it is immediate that 
$\stwomatr{\ell}{0}{0}{1} \Gamma(4p)  {\stwomatr{\ell}{0}{0}{1}}^{-1}
\subset \Gamma(4)$ if and only if $\ell \in \{1,p\}$.  Furthermore, the
space $\modforms_1(\Gamma(4))$ is stable under sending $f$ to $f| \gamma$
for $\gamma \in \Gamma(1)$.  Thus the only $\alpha$ that we need to
consider are representatives for the coset spaces
$\Gamma(1) \backslash \Gamma(1) \stwomatr{\ell}{0}{0}{1} \Gamma(1)$, for
$\ell \in \{1,p\}$.  Now a set of representatives for these is 
$S \union \left\{ \stwomatr{p}{0}{0}{1}\right\}$; hence our only remaining
task is to show that forms $f|\stwomatr{p}{0}{0}{1}$ do not contribute
further to the span.  This follows because the Hecke operator $T_p$ given
by the double coset $\Gamma(4) \stwomatr{p}{0}{0}{1} \Gamma(4)$ acts
linearly on $\modforms_1(\Gamma(4))$; hence $f|\stwomatr{p}{0}{0}{1}$ can
be written as a linear combination of $f|T_p$ (which, being in
$\modforms_1(\Gamma(4))$, is already in the span we are considering) and of
$f|\beta$ for other representatives
$\beta \notin \Gamma(4) \stwomatr{p}{0}{0}{1}$ of the coset space 
$\Gamma(4) \backslash \Gamma(4) \stwomatr{p}{0}{0}{1} \Gamma(4)$.  But all
these other representatives can be expressed as $\beta = \gamma \alpha$,
where $\alpha \in S$ and $\gamma \in \Gamma(1)$, and so the other terms
$f|\beta$ are already in the span that we are considering.
\end{proof}

\begin{remark}
\label{remark2.5}
The above lemma says that the maps from $X(4p)$ to $X(4)$ that are induced
by a map of the form $z \mapsto \alpha z$ on $\Half$ can all be obtained by
taking $\alpha  = \gamma\beta$, with 
$\beta \in S \union \{\stwomatr{p}{0}{0}{1}\}$, and $\gamma$ in a set
of representatives for $\Gamma(4) \backslash \Gamma(1)$.  The presence of
$\gamma$ corresponds to composing with an automorphism of $X(4)$.
\end{remark}

We now use the space $V$, viewed as an incomplete linear series on $X(4p)$,
to define a morphism $\morph$ from $X(4p)$ to projective space.
This morphism sends $z \in \Half$ to the point $\morph(z)$ with projective
coordinates coming from all the values
$\bigl( M|\alpha \bigr)(z)$,
$\bigl( N|\alpha \bigr)(z)$,
$\bigl( P|\alpha \bigr)(z)$,
for $\alpha \in S$, since these modular forms span $V$.  
We can write this concretely, but not entirely accurately, as
\begin{equation}
\label{equation2.4}
\morph(z) = \left[
M(z):N(z):P(z):
M|\stwomatr{1}{0}{0}{p}(z):
\dots
P|\stwomatr{1}{p-1}{0}{p}(z) \right].
\end{equation}
The problem with the above is that the spanning set 
$\bigl\{ M|\alpha, N|\alpha, P|\alpha \bigm| \alpha \in S \bigr\}$ of $V$
from Lemma~\ref{lemma2.4} is occasionally linearly dependent, a phenomenon
that we have observed experimentally whenever $p \equiv 3 \pmod{4}$.
However, this does not affect anything in our argument, and a purist should
use a basis for $V$ instead of the above spanning set.  At any rate,
we immediately see that $\morph(z)$ is invariant under the action of
$\Gamma(4p)$ on the upper half plane, and we can even take $z$ to belong to
the extended upper half plane $\Half^* = \Half \union \Projective^1(\Q)$,
so as to define $\morph$ on the cusps of $X(4p)$ as well.
Equation~\eqref{equation2.4} is incidentally well-defined because the first
three forms $M(z)$, $N(z)$, and $P(z)$ listed above already do not have any
common zeros on $\Half^*$, so that the incomplete linear series on $X(4p)$
associated to $V$ does not have base points.  Projecting the image of
$\morph$ onto these first three coordinates gives rise to the composition
of the map $X(4p) \to X(4)$, arising from $z \mapsto z$, with
the embedding $X(4) \to \Projective^2$ of Remark~\ref{remark2.2}.
Similarly, for every $\alpha \in S \union \stwomatr{p}{0}{0}{1}$, we can
project the image of $\morph$ to three coordinates corresponding to
$M|\alpha, N|\alpha, P|\alpha$, to obtain the composition of the map
$X(4p)\to X(4)$ arising from $z \mapsto \alpha z$ with the map $X(4)
\to \Projective^2$.  Thus our morphism $\morph$ contains at least all the
information that is to be found in the maps $X(4p) \to X(4)$ of
Remark~\ref{remark2.5}.  In fact, $\morph$ contains strictly more
information, as shown below when we match up factors of automorphy in the
proof of part~(2) of Theorem~\ref{theorem3.2}.

\section{Separation properties of the morphism $\morph$}
\label{section3}

We now discuss the extent to which $\morph$ fails to be an embedding of
$X(4p)$.  We first define a congruence group $\Gammapm(4p)$ that
contains $\Gamma(4p)$ with index~$2$:
\begin{equation}
\label{equation3.1}
\Gammapm(4p) = \left\{ \beta \in \Gamma(4) \mid \beta \equiv
\pm \twomatr{1}{0}{0}{1} \pmod{p} \right\}.
\end{equation}
Hence the corresponding modular curve $\Xpm(4p) = X(\Gammapm(4p))$ is the
quotient $X(4p)/\iota$ for an involution $\iota$, corresponding to the
action of a matrix $\iota \in \Gamma(4)$ such that $\iota \equiv
\stwomatr{-1}{0}{0}{-1} \pmod{p}$.  The following result is an immediate
consequence of the fact that for a primitive matrix $\alpha$ with $\det
\alpha = p$, we have $\alpha \Gammapm(4p) \alpha^{-1} \subset \Gamma(4)$:

\begin{lemma}
\label{lemma3.1}
The elements of $V$ all belong to the subspace $\modforms_1(\Gammapm(4p))
\subset \modforms_1(\Gamma(4p))$.  Hence the morphism $\morph$ factors
through $\Xpm(4p)$.
\end{lemma}

Hence, if two points of $X(4p)$ are related under $\iota$, then they cannot
be separated by the elements of $V$, which have proportionate values at any
pair of points $z, \iota z \in \Half$.  We can conversely show that this is
the only way that noncuspidal points of $X(4p)$ fail to be separated by
$\morph$.  In fact, viewing $\morph$ as a map from the quotient $\Xpm(4p)$
to projective space, the following theorem shows that $\morph$ is an
isomorphism away from the cusps of $\Xpm(4p)$.

\begin{theorem}
\label{theorem3.2}
\begin{enumerate}
\item
The elements of $V$ separate tangent vectors on $\Xpm(4p)$, i.e., for every
point $x \in \Xpm(4p)$, there exists an element of $V$ that vanishes at $x$
to order exactly~$1$.
\item
If $z, z' \in \Half$ satisfy $F(z) = F(z')$, then there exists
$\beta \in \Gammapm(4p)$ such that $z' = \beta z$.  Hence $z$ and $z'$
correspond to the same noncuspidal point of $\Xpm(4p)$, and the elements of
$V$ separate the noncuspidal points of $\Xpm(4p)$.
\end{enumerate}
\end{theorem}
\begin{proof}
\begin{enumerate}
\item
Observe that $\Gamma(4)$ does not contain any elliptic elements, and so the
map $\Xpm(4p) \to X(4)$ arising from $z \mapsto z$ is unramified except at
the cusps.  Hence, if $x \in \Xpm(4p)$ is not a cusp, we can already find
an element of $\mathrm{span}\{M,N,P\} = \modforms_1(\Gamma(4)) \subset V$
that vanishes at $x$ to order exactly~$1$, due to the projective embedding
$X(4) \to \Projective^2$ of Remark~\ref{remark2.2}.  On the other hand, if
$x$ is a cusp of $\Xpm(4p)$, then we may assume that $x = \infty$
because the subspace $V$ is stable under the action of $\Gamma(1) =
SL(2,\Z)$, which acts transitively on the cusps.  The stabilizer of
$\infty$ in the subgroup $\Gammapm(4p)$ is generated by
$\stwomatr{1}{4p}{0}{1}$, and so the form $N(z/p) \in V$ vanishes at the
cusp $\infty \in \Xpm(4p)$ to order precisely~$1$, since its
$q$-expansion begins with $4q^{1/4p}$.
\item
Since $F(z) = F(z')$ in projective coordinates, there must exist a nonzero
$\lambda \in \C^*$ such that for all $f \in \{M,N,P\}$ and $\alpha \in S$,
we have $f|\alpha(z') = \lambda f|\alpha(z)$.  (Actually, this holds for
all primitive $\alpha$ with $\det \alpha \in \{1,p\}$.)  We first
consider $\alpha = \stwomatr{1}{0}{0}{1}$.  Since
$(M(z'),N(z'),P(z')) = \lambda(M(z), N(z), P(z))$, we see from the
embedding $X(4) \to \Projective^2$ given by $[M:N:P]$ that that there
exists $\beta \in \Gamma(4)$ such that $z' = \beta z$; moreover,
$\lambda = j(\beta,z)$. 
At this point, however, we do not know whether
$\beta \in \Gammapm(4p)$.  Now repeat the above reasoning for arbitrary
$\alpha$.  Since $(M|\alpha(z'),N|\alpha(z'),P|\alpha(z'))
= \lambda(M|\alpha(z), N|\alpha(z), P|\alpha(z))$, we obtain that there
exists $\delta_\alpha \in \Gamma(4)$ for which
\begin{equation}
\label{equation3.2}
 \alpha z' = \delta_\alpha \alpha z, \quad \text{i.e., }
 \alpha \beta z = \delta_\alpha \alpha z.
\end{equation}
Comparing the actual values of the $f|\alpha$, and using the fact that
$f|\delta_\alpha = f$, we obtain that
$j(\delta_\alpha, \alpha z) j(\alpha, z')^{-1} 
   = \lambda j(\alpha, z)^{-1} = j(\beta, z) j(\alpha, z)^{-1}$.  
Hence
\begin{equation}
\label{equation3.3}
j(\delta_\alpha \alpha, z) 
= j(\delta_\alpha, \alpha z) j(\alpha, z)
= j(\beta, z) j(\alpha, z')
= j(\alpha \beta, z).
\end{equation}
Write $\delta_\alpha \alpha = \stwomatr{a}{b}{c}{d}$
and $\alpha \beta = \stwomatr{a'}{b'}{c'}{d'}$.  Then~\eqref{equation3.3}
implies that $c = c'$ and $d = d'$.
Similarly, combining with~\eqref{equation3.2}, we
obtain that $a = a'$ and that $b = b'$.  Thus $\delta_\alpha \alpha =
\alpha \beta$, and this shows that $\beta \in \alpha^{-1} \Gamma(4) \alpha$
for all $\alpha$.  It is now left to the reader to verify that
\begin{equation}
\label{equation3.4}
\Intersect_{\alpha \in S} \alpha^{-1} \Gamma(4) \alpha = \Gammapm(4p).
\end{equation}
A conceptual explanation of~\eqref{equation3.4} is that the identity
$\delta_\alpha \alpha = \alpha \beta$ for all primitive $\alpha$ with $\det
\alpha = p$ means that right multiplication by $\beta$ acts trivially on
the coset space $\Gamma(1)\backslash \Gamma(1) \stwomatr{p}{0}{0}{1}
\Gamma(1)$, which can be identified with the finite projective line
$\Projective^1(\F_p)$.  Thus $\beta$ must be congruent to a scalar matrix
modulo $p$, but $\det \beta = 1$, so the scalar matrix must be 
$\pm \stwomatr{1}{0}{0}{1}$.
\end{enumerate} 
\end{proof}

\begin{remark}
\label{remark3.3}
The above argument generalizes to the case of Shimura curves, in which case
we do not have to worry about cusps.  The analogous setting replaces
$\Gamma(4)$ by an arithmetic group $\Gamma^B$, coming from an indefinite
quaternion algebra $B$ over $\Q$.  We choose the level of $\Gamma^B$
appropriately to ensure that there are no elliptic elements, and
replace $\{M,N,P\}$ by a basis for a suitable space of modular forms
$\modforms_k(\Gamma^B)$ that yield a projective embedding of $X(\Gamma^B)$.
(This requires a sufficiently high weight~$k$.)  Analogs $\Gamma^B_p$ and
$\Gamma^B_{\pm p}$ of $\Gamma(4p)$ and $\Gammapm(4p)$, as well as of the
set $S$, exist for primes $p$ that do not divide either the level of
$\Gamma$ or the discriminant of $B$.  We can then generalize the statement
and proof of Theorem~\ref{theorem3.2} to a subspace 
$V^B \subset \modforms_k(\Gamma^B_{\pm p})$ that is analogous to the space
$V$.

\end{remark}

\begin{remark}
\label{remark3.3.5}
Generally, let $X$ be any variety over $\C$, and let $\morph$ be a morphism
from $X$ to some projective space $\Projective^N$.  Then it is well known
that if $\morph$ is an embedding, then the natural homomorphism
$H^0(\Projective^N, \mathcal{O}_{\Projective^N}(k))
       \to H^0(X, \morph^* \mathcal{O}_{\Projective^N}(k))$ 
is surjective for large $k$.
In our setting, with $X = \Xpm(4p)$ and our morphism $\morph$, the image of
the homomorphism is the space $W_k$ of Conjecture~\ref{conjecture1.1},
while $H^0(X, \morph^*\mathcal{O}_{\Projective^N}(k))
          = \modforms_k(\Gammapm(4p))$.
Now if it happened to be the case that $\morph$ were an embedding of
$\Xpm(4p)$, then we would be able to conclude that for large $k$, 
$W_k = \modforms_k(\Gamma_\pm(4p))$.  We also know the dimension of the
space of modular forms of weight $k \geq 2$:
\begin{equation}
\label{equation3.4.5}
\dim \modforms_k(\Gammapm(4p)) = (1/2) \dim \modforms_k(\Gamma(4p))
= (p^2 - 1)((k-1)p + 3/2)
\end{equation}
because the degree two map from $X(4p)$ to $\Xpm(4p)$ is unramified
everywhere, including the cusps.  Thus we would obtain the result of
Conjecture~\ref{conjecture1.1} for large $k$, but without the correction
term $3(p+1)(p-3)$.  We shall however see that $\morph$ partially
identifies the cusps of $\Xpm(4p)$, making the correction term necessary if
the conjecture is to have any chance of being true.  However, in the case
of quaternion algebras, we see by Remark~\ref{remark3.3} that the analog of
Conjecture~\ref{conjecture1.1} is indeed true for large $k$, without a
correction term.  As for the question of how large $k$ must be, our
conjecture that $k \geq 4$ is sufficient is considerably smaller than the
bound proved in~\cite{GrusonLazarsfeldPeskine} for arbitrary curves.  Our
hope is that the behavior of our linear system $V$ on the modular curve
$\Xpm(4p)$ is considerably better than the worst case.  Compare this to the
result of~\cite{BorisovGunnellsHigherWeight}, which states that the cuspidal
parts of the weight~$k$ monomials in all Eisenstein series of weight~1 on
$\Gamma_1(N)$ generate the full space of cusp forms
$\mathcal{S}_k(\Gamma_1(N))$, for $k \geq 3$ and $N \geq 5$. 
\end{remark}

To describe the effect of our morphism $\morph$ on the cusps of
$\Xpm(4p)$, we first introduce some notation.  Consider a cusp
$a/c \in \Projective^1(\Q)$ in the extended upper half plane, where
$\gcd(a,c)=1$.  We choose $b,d \in \Z$ such that 
$\gamma = \stwomatr{a}{b}{c}{d}$ belongs to $\Gamma(1)$, so our cusp is
$a/c = \gamma \infty$.  We usually represent our cusp by the primitive
vector $\stwovec{a}{c} = \gamma \stwovec{1}{0} \in \Z^2$.  Our cusp can
also be represented by $\stwovec{-a}{-c}$; we shall point out the occasions
when we may need to change signs, but this ambiguity is generally harmless.
We define the following equivalence relations on primitive vectors, which
also give equivalence relations on the set $\Projective^1(\Q)$ of cusps:
\begin{equation}
\label{equation3.5}
\begin{split}
\twovec{a'}{c'} \sim \twovec{a}{c} &\iff
     \exists \mu \in (\Z/4p\Z)^* \text{ such that } 
        \twovec{a'}{c'} \equiv \mu \twovec{a}{c} \pmod{4p}, \\
\twovec{a'}{c'} \approx \twovec{a}{c} &\iff
        \twovec{a'}{c'} \equiv \pm \twovec{a}{c} \pmod{4} \text{ and }
        \twovec{a'}{c'} \equiv \pm \twovec{a}{c} \pmod{p}.
\end{split}
\end{equation}
Thus the finer relation $\approx$ restricts us to the four values of
$\mu \in (\Z/4p\Z)^*$ satisfying $\mu \equiv \pm 1 \pmod{p}$, since in all
cases we have $\mu \equiv \pm 1 \pmod{4}$.  We can easily count the number
of equivalence classes of cusps under each relation:
\begin{equation}
\label{equation3.6}
\abs{\Projective^1(\Q)/\sim} = 6(p+1),
\qquad
\abs{\Projective^1(\Q)/\approx} = 3(p^2 - 1).
\end{equation}

\begin{lemma}
\label{lemma3.4}
\begin{enumerate}
\item
Two cusps $\stwovec{a}{c}$ and $\stwovec{a'}{c'}$ map to the same element
of $\Xpm(4p)$ if and only if $\stwovec{a'}{c'} \approx \stwovec{a}{c}$.
\item
Two cusps $\stwovec{a}{c}$ and $\stwovec{a'}{c'}$ map to the same element
of $X(4)$ if and only if
$\stwovec{a'}{c'} \equiv \pm \stwovec{a}{c} \pmod{4}$.
\end{enumerate}
\end{lemma}
\begin{proof}
This is standard; see for example Section 1.6 of \cite{Shimura}.
\end{proof}

Even though cusps on the curve $\Xpm(4p)$ are determined up to equivalence
classes under $\approx$, our next result states that the morphism
$\morph$ only distinguishes cusps up to the coarser equivalence relation
$\sim$.

\begin{theorem}
\label{theorem3.5}
Two cusps $\stwovec{a}{c}$ and $\stwovec{a'}{c'}$ have the same image under
$\morph$ if and only if $\stwovec{a}{c} \sim \stwovec{a'}{c'}$.
\end{theorem}
\begin{proof}
We first show the ``if'' direction.  Assume that 
$\stwovec{a}{c} \sim \stwovec{a'}{c'}$, where 
$\stwovec{a}{c} = \gamma \stwovec{1}{0}$ for some $\gamma \in \Gamma(1)$.
Possibly after a sign change, we can assume that 
$\stwovec{a'}{c'} \equiv \stwovec{a}{c} \pmod{4}$.  It follows
(this is essentially part~2 of the previous lemma)
that there exists $\beta \in \Gamma(4)$ satisfying
$\stwovec{a'}{c'} = \beta \stwovec{a}{c} = \beta \gamma \stwovec{1}{0}$.
Now in order to show that the two cusps have the same image under $\morph$,
we must show that there exists $\lambda \in \C^*$ such that for all forms
$f|\alpha$, with $f \in \{M,N,P\}$ and $\alpha \in S$, we have that the
constant term in the $q$-expansion of $f|\alpha\beta\gamma$ is $\lambda$
times the constant term of $f|\alpha\gamma$.  (Taking the constant terms
amounts to ``evaluating'' the form $f|\alpha$ at the cusps 
$a'/c' = \beta\gamma\infty$ and $a/c = \gamma\infty$.)  From the
choice $\alpha = \stwomatr{1}{0}{0}{1}$ and the fact that $f|\beta = f$,
we obtain that $\lambda$ must be equal to $1$.  Now consider 
$\alpha = \stwomatr{1}{t}{0}{p}$ for $t \in \{0, \dots, p-1\}$.  Write
$\alpha\gamma\infty = (a+tc)/pc = A/C$ in lowest terms; in other words, 
$\alpha\gamma\stwovec{1}{0} = \alpha\stwovec{a}{c} = \stwovec{a+tc}{pc} 
= g \stwovec{A}{C}$, where $g = \gcd(a+tc, pc) \in \{1,p\}$.  Similarly
let $\alpha\beta\gamma\infty = A'/C'$, with $g' \in \{1,p\}$ and
$\alpha\beta\gamma\stwovec{1}{0} = g' \stwovec{A'}{C'}$.  
Now observe that $g = p$ if and only if $a+tc \equiv 0 \pmod{p}$, which is
equivalent (since $\stwovec{a}{c} \sim \stwovec{a'}{c'}$)
to $a'+tc' \equiv 0 \pmod{p}$, which happens if and only
if $g' = p$; thus $g = g'$.  As a first consequence,
$\stwovec{A}{C} \equiv \stwovec{A'}{C'} \pmod{4}$, which implies that there
exists $\delta_\alpha \in \Gamma(4)$ such that
$\stwovec{A'}{C'} = \delta_\alpha \stwovec{A}{C}$.  Multiply this equation
by $g=g'$ to obtain
\begin{equation}
\label{equation3.7}
\alpha\beta\gamma \twovec{1}{0} = \twovec{a'+tc'}{pc'} 
 = \delta_\alpha \twovec{a+tc}{pc}
 = \delta_\alpha \alpha \gamma \twovec{1}{0}.
\end{equation}
Taking into account the fact that $\det \alpha\beta\gamma = p
= \det \delta_\alpha \alpha\gamma$, we see that
\begin{equation}
\label{equation3.8}
\exists t_\alpha \in \Q \text{ such that } 
  \alpha \beta \gamma
       = \delta_\alpha \alpha \gamma \twomatr{1}{t_\alpha}{0}{1}.
\end{equation}
Thus for $f\in\{M,N,P\}$, we obtain that $f|\alpha\beta\gamma =
f|\alpha\gamma\stwomatr{1}{t_\alpha}{0}{1}$ has the same constant term as
$f|\alpha\gamma$, as desired.  

We now show the converse (``only if'').  Let $a/c$ and $a'/c'$ be cusps
with the same image under $\morph$.  As in the proof of
Theorem~\ref{theorem3.2}, the cusps must be related by an element $\beta
\in \Gamma(4)$.  Take $\gamma \in \Gamma(1)$ such that $\stwovec{a}{c}
= \gamma\stwovec{1}{0}$; possibly after a harmless sign change, we can
assume then that $\stwovec{a'}{c'} = \beta \stwovec{a}{c}
= \beta\gamma \stwovec{1}{0}$.  Thus, by looking as before at 
$f \in \{M,N,P\}$ and $\alpha = \stwomatr{1}{0}{0}{1}$, we obtain that
$\lambda = 1$, where $\lambda$ is the ratio between the constant terms of
$f|\alpha\beta\gamma$ and $f|\alpha\gamma$.
Now repeat the same argument that led 
us to~\eqref{equation3.2} (applied to the cusp $a/c = \gamma\infty$ instead
of to $z \in \Half$) to obtain that for each $\alpha \in S$, there
exists $\delta_\alpha \in \Gamma(4)$ for which $\alpha\beta\gamma\infty =
\delta_\alpha \alpha\gamma\infty$.  Taking the stabilizer of $\infty$ into
account, we obtain an identity that is at first glance slightly weaker
than~\eqref{equation3.8}: 
\begin{equation}
\label{equation3.9}
\exists r_\alpha, s_\alpha, t_\alpha \in \Q \text{ such that } 
  \alpha \beta \gamma
   = \delta_\alpha \alpha \gamma \twomatr{r_\alpha}{t_\alpha}{0}{s_\alpha}.
\end{equation}
Comparing constant terms, this implies that $s_\alpha = \lambda^{-1} = 1$.
Comparing determinants, we also conclude that $r_\alpha = 1$, so we do
obtain~\eqref{equation3.8} after all.  We can now run the argument from the
first part of this proof backwards: first conclude~\eqref{equation3.7}, and
deduce that for all $t \in \{0, \dots, p-1\}$ we have the equality
$\gcd(a'+tc', pc') = \gcd(a+tc, pc)$.  Hence there exists $\mu$ for which
$\stwovec{a'}{c'} \equiv \mu \stwovec{a}{c} \pmod{p}$.  Since 
$\beta \equiv \stwomatr{1}{0}{0}{1} \pmod{4}$, we 
obtain in addition that $\stwovec{a'}{c'} \equiv \stwovec{a}{c} \pmod{4}$.
This shows that $\stwovec{a'}{c'} \sim \stwovec{a}{c}$, as desired.
\end{proof}

\begin{corollary}
\label{corollary3.6}
For all $k \geq 3$, the dimension of the space $W_k$ of
Conjecture~\ref{conjecture1.1} is less than or equal to
$(1/2) \dim \modforms_k(\Gamma(4p)) - 3(p+1)(p-3)$.
\end{corollary}
\begin{proof}
By the above theorem, elements of $V$, and hence of $W_k$, always have the
same value (in the sense of having the same constant term) at all cusps
related by $\sim$.  However, there are enough Eisenstein series in the full
space $\modforms_k(\Gammapm(4p))$ to know that modular forms in the full
space can have arbitrary values at the cusps.
By~\eqref{equation3.6}, this means that when we restrict our attention to
forms in $W_k$, the relations between the values at the cusps impose
$3(p^2-1) - 6(p+1) = 3(p+1)(p-3)$ independent linear conditions --- hence
the correction term in Conjecture~\ref{conjecture1.1}.
\end{proof}

\section{Empirical results}
\label{section4}

In this section we report on the computer calculations that led us to
Conjecture~\ref{conjecture1.1}.  Let $\{f_1, \dots, f_d\}$ be a basis for
our space $V$, with $d \leq 3(p+1)$.  Then the space $W_k$ of
Conjecture~\ref{conjecture1.1} can be viewed as the degree~$k$ part of the
image of the morphism of graded rings
\begin{equation}
\label{equation4.1}
\varphi: \C[x_1, \dots, x_d]
             \to \Directsum_{k \geq 0} \modforms_k(\Gammapm(4p)),
\qquad\qquad\text{with } \varphi(x_i) = f_i.
\end{equation}
We shall represent a given form $g(z) \in \modforms_k(\Gammapm(4p))$ by its
$q$-expansion in terms of $q^{1/4p} = \exp(2\pi i z/4p)$.  We only need to 
consider a finite number of coefficients, and can choose a number $L$
(depending on $p$ and $k$) large enough so that we have an injective linear
transformation
\begin{equation}
\label{equation4.2}
g(z) = \sum_{n=0}^\infty a_n q^{n/4p} \mapsto
   (a_0, \dots, a_{L-1}) \in \C^L.
\end{equation}
Thus finding the dimensions of the $W_k$ reduces to a linear algebra
calculation in $\C^L$.  Our first observation is that the calculation
can be carried out entirely over $\Q$, or even over $\Z$, by clearing
denominators.  Moreover, the coefficients $a_n$ above can be grouped
according to the residue class of $n \bmod 4p$.

\begin{lemma}
\label{lemma4.1}
We can choose a basis $\{f_1, \dots, f_d\}$ for $V$ that consists of forms
$f_i = \sum_{n=0}^\infty a_n q^{n/4p}$ satisfying the following two
conditions: 
\begin{enumerate}
\item
all the coefficients satisfy $a_n \in \Q$,
\item
there exists $b \in \Z/4p\Z$ such that $a_n \neq 0$ only if $n \equiv b
\pmod{4p}$.
\end{enumerate}
Consequently, the monomials $f_1^{m_1} \dots f_d^{m_d}$ of total degree
$k = \sum_{i} m_i$ that span $W_k$ will all have $q$-expansions satisfying
properties (1) and~(2) above.
\end{lemma}
\begin{proof}
It is sufficient to find a set of generators of $V$ satisfying the above
two conditions.  By~\eqref{equation2.1}, every form $f \in \{M,N,P\}$ can
be expanded as $f = \sum_n c_n q^{n/4}$, with $c_n \in \Z$, such that $c_n
\neq 0$ only if $n$ belongs to a fixed residue class mod~$4$.  Upon viewing
this $q$-expansion in terms of $q^{1/4p}$, we obtain that $M$, $N$, and $P$
all satisfy the above conditions.  Now the rest of $V$ is spanned by the
forms $f^{(4t)} = f((z+4t)/p)$ for $f \in \{M,N,P\}$ and $t \in \Z/p\Z$,
but these forms have $q$-expansions involving $p$th roots of unity.  (It is 
technically useful here to use $4t$ instead of $t$, which amounts to using
a slightly different set of coset representatives for 
$\Gamma(1) \backslash \Gamma(1) \stwomatr{p}{0}{0}{1} \Gamma(1)$
than the representatives that we used for our original set $S$.  The forms
$f((z+t)/p)$ would have involved working with $(4p)$th roots of unity.)
We can replace the $f^{(4t)}$ with the collection of modular
forms $f_{(b)}$, for $b \in \Z/p\Z$, given by
\begin{equation}
\label{equation4.3}
f_{(b)} = \sum_{n \equiv b \pmod{p}} c_n q^{n/4p} =
(1/p)\sum_{t \in \Z/p\Z} f^{(4t)} \exp(-2\pi i bt/p).
\end{equation}
Thus $f_{(b)}$ selects only the terms in the
$q$-expansion of $f(z/p)$ for which $n$ belongs to a fixed residue class
mod~$p$.  Since we also have a congruence condition on $n \bmod 4$, we see 
that the $f_{(b)}$ satisfy both of the two conditions above.
\end{proof}

The above lemma not only allows us to do our entire calculation over $\Z$, 
but also allows us to decompose each vector space
$\modforms_k(\Gammapm(4p))$ and $W_k$ into a direct sum of $4p$  
subspaces, depending on the parameter $b \in \Z/4p\Z$.  This means
that instead of doing linear algebra on a matrix of size $L$ (i.e., with
approximately $L$ rows and columns), we can instead work with $4p$ smaller
matrices, each of which has size $L/4p$.  This is more efficient: for
example, if we were to use brute-force Gaussian elimination without
noticing the zero entries in our matrix, this would reduce the number of
arithmetic operations from $O(L^3)$ to $O(4p(L/4p)^3)$, a factor of
$16p^2$.  Conceptually, our approach amounts to decomposing our spaces
of modular forms into isotypic components under the action of the cyclic
group $\{\stwomatr{1}{b}{0}{1}\mid b \in \Z/4p\Z\}$, acting on modular
forms.  It would be interesting to extend this approach by taking into 
account the full action of the finite group
$PSL(2,\Z/4p\Z) \isomorphic \Gamma(1)/\Gammapm(4p)$.

We can now describe how we proceeded to find bases for the spaces $W_k$ for
$k \in \{1, 2, 3\}$, and, when practical, for $k=4$.  We directly computed
linear dependencies among the $3(p+1)$ forms that generate $V$, to obtain a
basis $\{f_1, \dots, f_d\}$ for $W_1 = V$.  We also directly computed
dependencies among the monomials in the $f_i$ of degree~$2$, in order to
find a basis for $W_2$.  The dependencies that we found were thus
(generators for) the degree~$2$ elements of the homogeneous ideal 
$I = \ker \varphi$, in terms of the homomorphism $\varphi$
of~\eqref{equation4.1}.  In order to compute $W_3$, we did not simply list
all the degree~$3$ monomials in the $f_i$; rather, we took into account the
ideal $I_2 \subset I$ generated by the degree~$2$ elements of $I$ that we
had already encountered.  We computed a Gr\"obner basis of $I_2$, which
contained both the above mentioned degree~$2$ elements of $I_2$ and some
elements in degree~$3$, and this allowed us to eliminate from consideration
several monomials which we already knew would be redundant in $W_3$, before
doing the linear algebra.
Once we found a basis for $W_3$, we obtained again generators of $I_3$, the
subideal of $I$ generated by the elements of degree $\leq 3$.  This time,
it was not practical to compute a Gr\"obner basis for $I_3$, but we were
still able to use our existing elements of degree $\leq 3$ to eliminate
many degree~$4$ monomials before doing the linear algebra needed to compute
the dimension of $W_4$.

\begin{table}
\begin{tabular}{|c|c|c|c|c|}
\hline
 $p$&$k=1$& $k=2$       &  $k=3$        &  $k=4$              \\
\hline
\hline
  3 &   9 &   33   (36) &    60    (60) &    84           (84)\\  
\hline
  5 &  18 &  115  (120) &   240   (240) &   360          (360)\\  
\hline
  7 &  21 &  189  (312) &   630   (648) &   984          (984)\\  
\hline
 11 &  33 &  499 (1212) &  2532  (2532) &  3852         (3852)\\  
\hline
 13 &  42 &  842 (2016) &  4200  (4200) &  6384         (6384)\\  
\hline 
 17 &  54 & 1359 (4572) &  9333  (9468) & ?????        (14364)\\  
\hline
 19 &  57 & 1624 (6420) & 13260 (13260) &\quad$\star$\quad\  (20100)\\  
\hline
   23 &   69 &   2347   (11496) & 23442  (23640) & ?????         (35784)\\  
\hline
\end{tabular}
\par\smallskip\noindent
\textbf{Table 1:}
This shows our computations of $\dim W_k$ for small $p$ and $k$.
\par\noindent
The numbers in parentheses are the upper bounds from
Corollary~\ref{corollary3.6}.
\end{table}

The results of our computations are shown in Table 1. 
The table lists the dimensions of $W_k$ that we were able to calculate for
each $p$, and also mentions in parentheses the upper bound on the dimension
from Corollary~\ref{corollary3.6}, namely 
$(k-1)p^3 - 3p^2/2 + (7-k)p + 15/2$, that we hope by
Conjecture~\ref{conjecture1.1} is equal to 
$\dim W_k$ for $k \geq 4$.  Note that $W_k$ already ``fills out'' the
conjectured dimension in weight~$3$ for all the primes except for $7$,
$17$, and~$23$ in our table.  Even for those primes, the deficit for $k=3$
is small, which leads us to be optimistic that the full dimension is
attained in weight~$4$.
We put a star ($\star$) for $p=19, k=4$ instead of question marks because
we are confident that the dimension is correct there, our bound
having already been attained in weight~$3$ for $p=19$.  We finally remark
that there is no possibility of attaining our bound on the dimension for 
$k \leq 2$, since $\dim W_1 \leq 3(p+1)$, and $\dim W_2$ grows at most
quadratically in $p$ since it is bounded above by the number of monomials
of degree~$2$ in the $f_i$.


\bibliographystyle{amsalpha}


\providecommand{\bysame}{\leavevmode\hbox to3em{\hrulefill}\thinspace}
\providecommand{\MR}{\relax\ifhmode\unskip\space\fi MR }
\providecommand{\MRhref}[2]{%
  \href{http://www.ams.org/mathscinet-getitem?mr=#1}{#2}
}
\providecommand{\href}[2]{#2}


\end{document}